\documentclass{article}
\usepackage{amsfonts}
\usepackage{amsthm}
\usepackage{latexsym}
\usepackage{amssymb}
\usepackage{amstext}
\usepackage{amsmath}

\usepackage{graphicx}
\usepackage{amssymb}
\usepackage{epstopdf}
\theoremstyle{plain}

\newtheorem{prop}{Proposition}

\theoremstyle{definition}
\newtheorem{defn}{Definition}

\def\G{\overline{G}}

\begin{document}

\begin{center}
\LARGE {\bf \textsc{The Homomorphism Poset of $K_{3,3}$} }
\end{center}\bigskip

\begin{center}
\textsc{Sally Cockburn}

\textsc{Department of Mathematics}

\textsc{Hamilton College,  Clinton,  NY 13323}

\textsc {\sl scockbur@hamilton.edu}

\end{center}


\begin{abstract}
 A {\it geometric graph} $\G$ is a simple graph drawn in the plane,  on points in general position, with straight-line edges.   We call $\G$  a {\it geometric realization} of the underlying abstract graph $G$.  A {\it geometric homomorphism} $f:\overline{G} \to \overline{H}$ is a vertex map that preserves adjacencies and crossings (but not necessarily non-adjacencies or non-crossings). Geometric homomorphisms can be used to define a partial order on the set of isomorphism classes of geometric realizations of an abstract graph $G$.  In this paper, the homomorphism poset of $K_{3,3}$ is determined.
 \end{abstract}

\section{Geometric Realizations} 

In \cite{HH}, Harborth defines a {\it good} drawing of a graph $G$ to be a drawing in the plane in which any two edges intersect at most once, and no three edges intersect at a common point; however, edges need not be represented with straight lines.  He further defines any two such drawings to be isomorphic if and only if there exists a graph isomorphism that preserves edge crossings and non-crossings, as well as regions and parts of edges.  By these definitions, there are 102 non-isomorphic good drawings of $K_{3,3}$.  Harborth also proves in this paper that if $m \equiv n \equiv 1 \mod 2$, then the parity of the number of crossings in any good drawing of $K_{m,n}$ is the same. 

A {\it geometric realization} (or {\it  rectilinear drawing}) of a graph $G$ is a drawing in the plane in which vertices are in general position and all edges are represented by straight lines.  Two realizations of $G$ are isomorphic if and only if there exists a graph isomorphism that preserves edge crossings and non-crossings; this is  weaker concept than Harborth's because it does not take into account regions and parts of edges.

\begin{figure}[htbp] 
   \centering
   \includegraphics[width=5.5in]{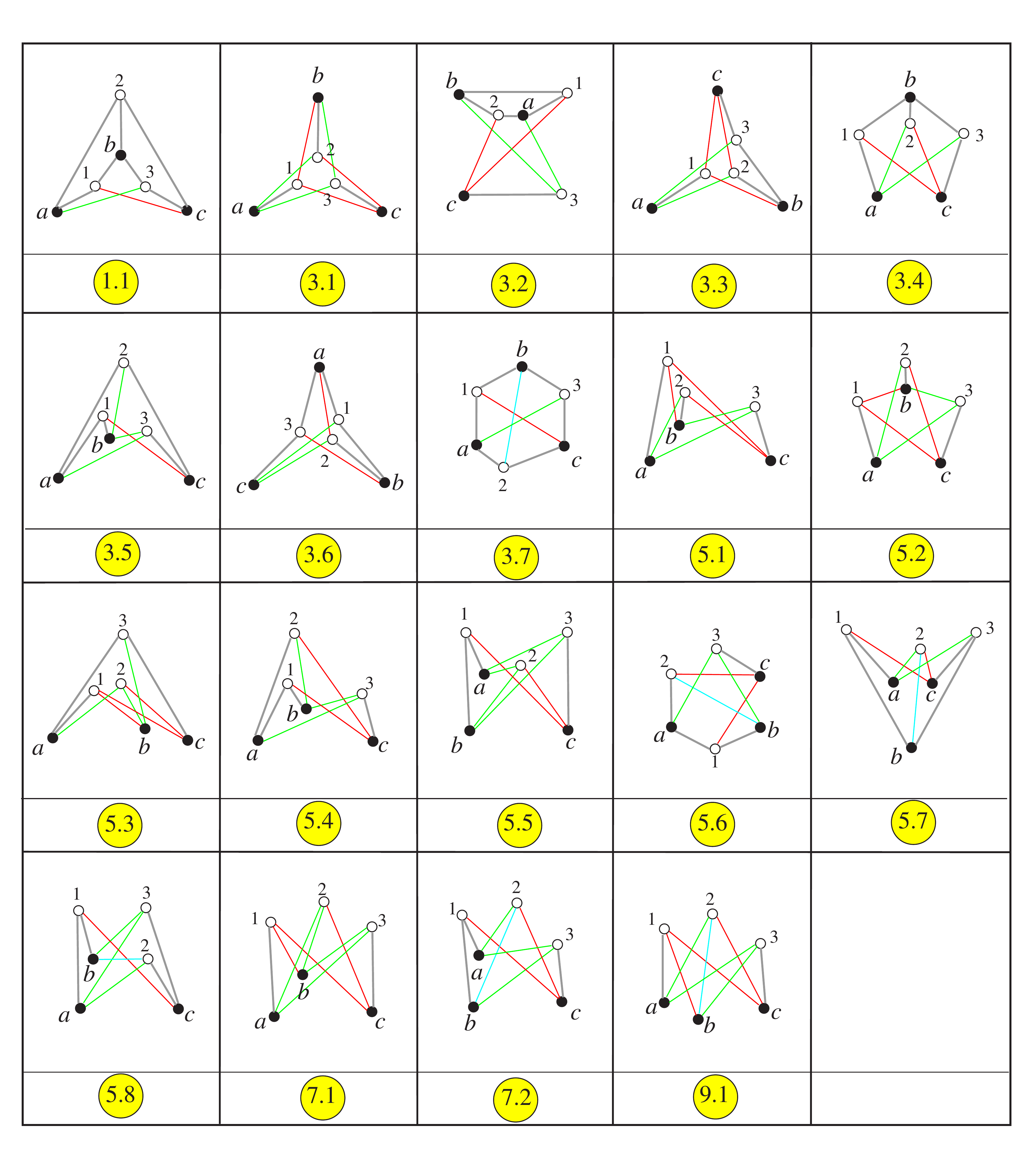} 
   \caption{The nineteen geometric realizations of $K_{3,3}$.}
   \label{fig:allK33}
\end{figure}

To determine the number of non-isomorphic geometric realizations of $K_{3,3}$,  note that any geometric realization of $K_{3,3}$ can be completed to obtain a geometric realization of $K_6$.  
As shown in \cite{BCDM}, there are 15 different geometric realizations of $K_6$, and for each of these, there are 10 ways of dividing the labeled vertices into two partite sets, but not all of these will result in different  realizations of $K_{3,3}$. The 19 different geometric realizations of $K_{3,3}$ are given in Figure~\ref{fig:allK33}.  Since any geometric realization is a good drawing, Harborth's  result also explains why the number of crossings is always odd.
\bigskip

To demonstrate that these realizations are non-isomorphic, we make use of some results from \cite{BCDM}.  First we recall some geometric graph invariants.

\begin{defn} Let $\overline{G}$ be a geometric realization of a  graph $G$. Then 
\begin{enumerate}
\item $cr(\G)$ is the total number of edge crossings in $\G$;  
\item for all $e\in E(\G)$,  $cr(e)$ is the number of edges crossing  $e$ in $\G$;
\item  $E_0 = \{e \in E(\overline{G}) \mid cr(e)=0\}$ is the set of uncrossed edges; 
\item the {\it uncrossed subgraph} of $\G$ is the abstract graph $\G_0 = (V(G), E_0)$;   
\item the {\it edge crossing graph} of $\G$  is the abstract graph $EX(\G)$ whose vertices are the edges of $G$, with adjacency when the corresponding edges of $\G$ cross; 
\item the {\it line/crossing graph} of $\G$ is  the $2$-edge colored abstract graph  $LEX(\G)$  whose vertices are the edges of $G$, with solid edges corresponding to the edges of $EX(\G)$ and dashed edges corresponding to the edges of the line graph, $L(G)$ (indicating when two edges of $G$ are adjacent).
\end{enumerate}
\end{defn}

\begin{figure}[htbp] 
   \centering
   \includegraphics[width=5.5in]{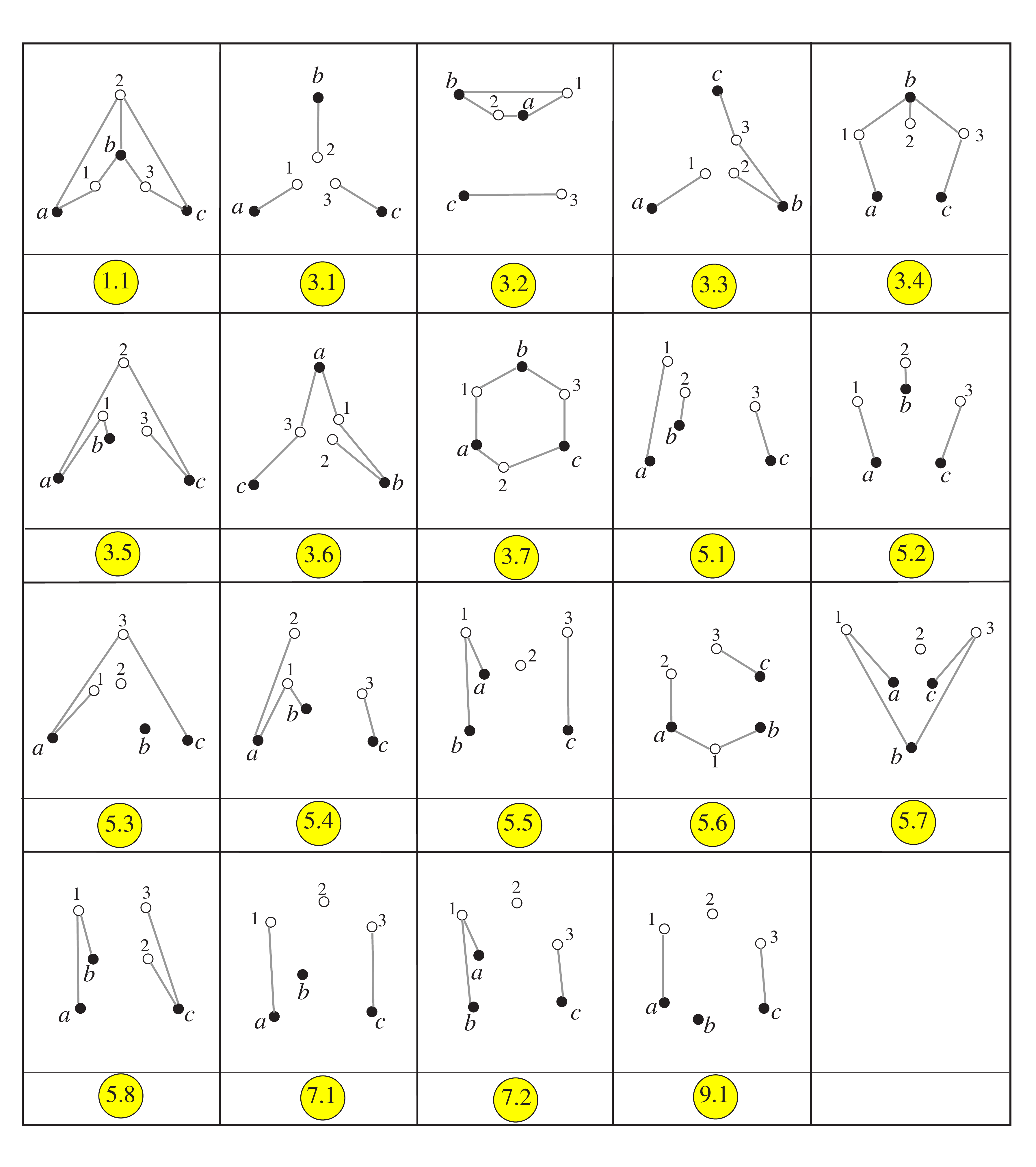} 
   \caption{The uncrossed subgraphs of the realizations in Figure~\ref{fig:allK33}.}
   \label{fig:allK330}
\end{figure}

\begin{figure}[htbp] 
   \centering
   \includegraphics[width=5.5in]{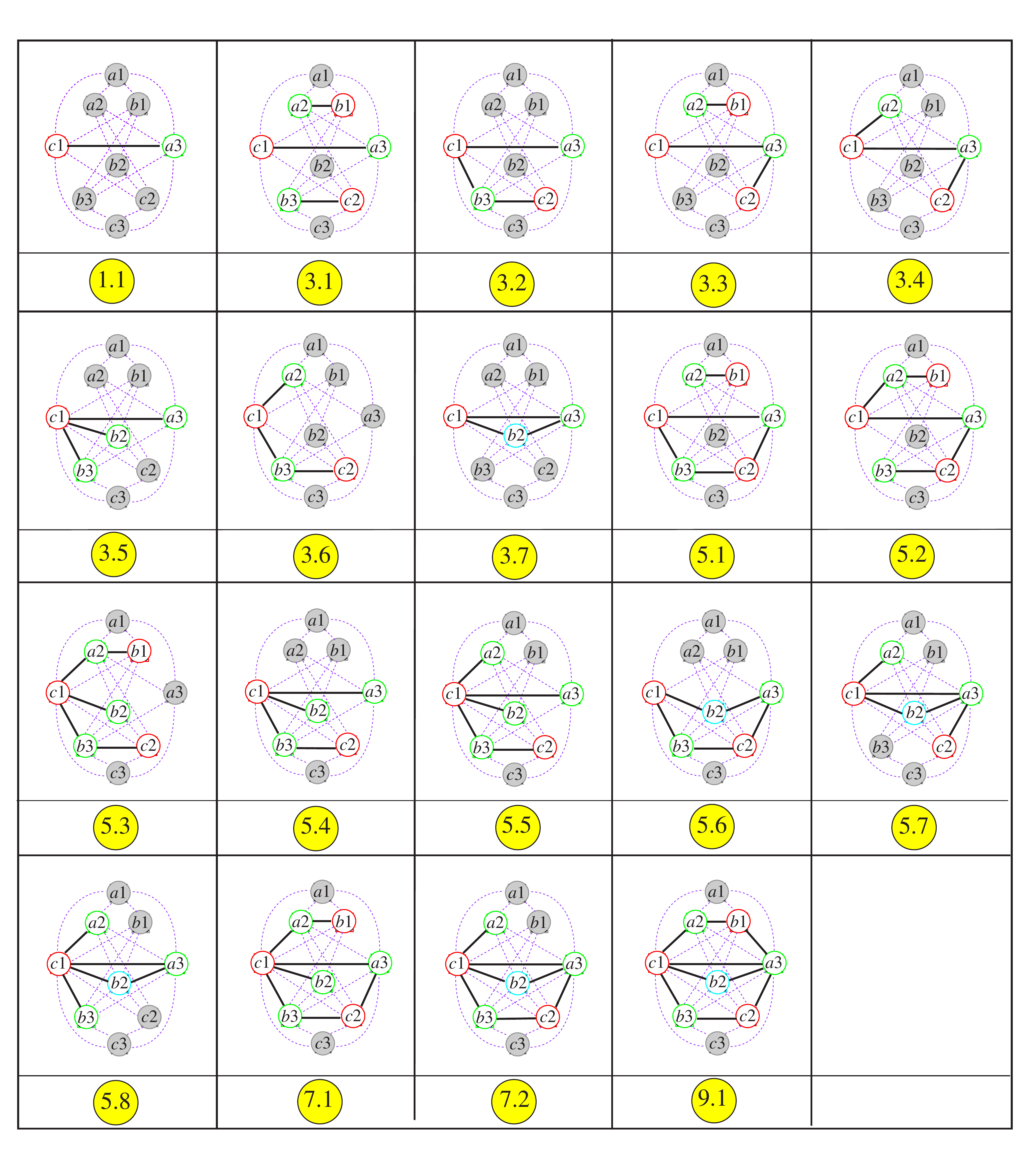} 
   \caption{The line/crossing graphs of the realizations in Figure~\ref{fig:allK33}.}
   \label{fig:lineXK33}
\end{figure}

Figure~\ref{fig:allK330}   gives the uncrossed subgraphs  and Figure~\ref{fig:lineXK33}  the line/crossing graphs of the realizations in Figure~\ref{fig:allK33}.

It suffices to show that realizations with the same number of edge crossings are not isomorphic.
From Figure~\ref{fig:allK330}, note that all realizations of $K_{3,3}$ with 3 crossings have different uncrossed subgraphs, except for realizations 3.5 and 3.6, both of which have uncrossed subgraph  $P_6$.  However, the edge crossing graph of realization 3.5 has a vertex of degree 3, while that of 3.6 has maximum degree 2.

Moving on to realizations of $K_{3,3}$ with 5 crossings,  realizations 5.1 and 5.2 both have $3K_2$ as uncrossed subgraph, but 5.1 has edge crossing graph $C_4 \cup K_2 \cup 3K_1$, whereas 5.2 has edge crossing graph $P_6 \cup 3K_1$. Realizations 5.4 and 5.6 both have $P_4 \cup K_2$ as uncrossed subgraph, but the edge crossing graph of the former is a tree and of the latter is a 5-cycle.

\section{Poset Structure}

Geometric homomorphisms were introduced in \cite{BC} as a natural generalization of abstract graph homomorphisms. A {\it geometric homomorphism} $f:\overline{G} \to \overline{H}$ is a vertex map that preserves adjacencies and crossings, but not necessarily non-adjacencies or non-crossings.  If $\G$ and $\widehat{G}$ are geometric realizations of $G$,  set $\G\preceq \widehat G$ if and only if there is a vertex-injective geometric homomorphism $f:\G \to \widehat G$. It is easy to verify that this defines a partial order on the set of all isomorphism classes of geometric realizations of $G$; the resulting structure is called the {it homomorphism poset} of $G$, denoted $\mathcal{G}$.  

\begin{figure}[htbp] 
   \centering
   \includegraphics[width=4in]{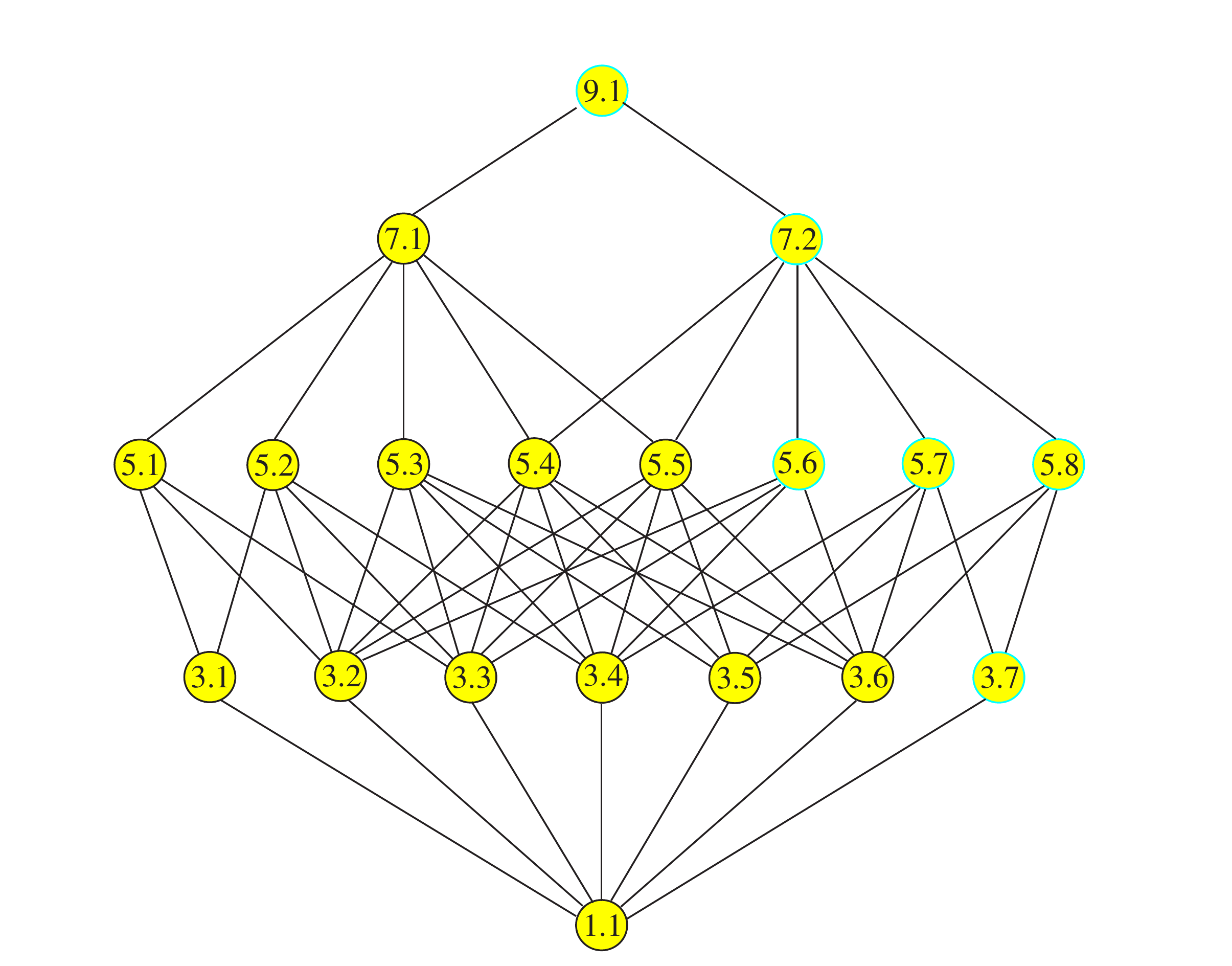} 
   \caption{The Hasse diagram of $\mathcal{K}_{3,3}$.}
   \label{fig:HasseK33}
\end{figure}

The Hasse diagram for the poset $\mathcal{K}_{3,3}$ is given in Figure~\ref{fig:HasseK33}. The nodes with  blue circumferences correspond to realizations of edge thickness 3. 

In \cite{BCDM}, it is shown that the homomorphism posets of other graphs of order $n=6$, namely $P_6, C_6$ and $K_6$, are neither lattices nor graded posets. However, it is obvious from Figure~\ref{fig:HasseK33} that $\mathcal{K}_{3,3}$ is a graded poset with rank function $\rho(\G) = \lfloor cr(\G)/2 \rfloor$. However, $\mathcal{K}_{3,3}$ is not a lattice; for example, realizations 3.1 and 3.2 have both realizations 5.1 and 5.2 as suprema.
Since the poset has a unique maximum, every geometric graph that is homomorphic to some realization of $K_{3,3}$ is homomorphic to realization 9.1; note also that every geometric graph that is homomorphic to a realization of $K_{3,3}$ of edge thickness 2 is homomorphic to realization 7.1.

Most edges in the Hasse diagram of $\mathcal{K}_{3,3}$ are induced by the identity map. The only ones  that aren't are from realizations with 3 edge crossings to those with 5 edge crossings.  Details are given in  Table~\ref{table:homdefs} and Table~\ref{table:homs}.

\begin{table}[htdp]
\caption{Non-identity isomorphisms on $K_{3,3}$.}
\begin{center}
\begin{tabular}{|c|ccccccc|} 
\hline
    & $f_1$ &  $f_2$ & $f_3$ & $f_4$ & $f_5$ & $f_6$ & $f_7$ \\ \hline
$a$ & $b$  & $c$ & $c$ & $b$ & $a$ & $b$ & 3 \\
$b$ & $a$  & $a$ & $b$ & $a$ & $b$ & $a$ & 1 \\
$c$ & $c$  & $b$ & $a$ & $c$ & $c$ & $c$ & 2 \\
1 & 2  & 3 & 3 & 1 & 1 & 1 & $b$ \\
2 & 1  & 2 & 2 & 2 & 3 & 3 & $a$ \\
3 & 3  & 1 & 1 & 3 & 2 & 2 & $c$  \\ \hline
\end{tabular}
\end{center}
\label{table:homdefs}
\end{table}%

\begin{table}[htdp]
\caption{Geometric homomorphisms from  3-crossing to 5-crossing realizations.}
\begin{center}
\begin{tabular}{|c|cccccccc|}
\hline
realization& 5.1 & 5.2 & 5.3 & 5.4 & 5.5 & 5.6 & 5.7 & 5.8 \\ 
\hline
3.1 & id & id &  &  &  &  &   & \\
3.2 & id & $f_1$ & $f_1$ & id & id & id & & \\
3.3 & id & id & $f_1$ & $f_2$ & $f_3$ & $f_2$ & & \\
3.4 & & id & $f_4$ & $f_4$ & $f_4$ & $f_4$ & id & \\
3.5 &  & & $f_5$ & id & id &  & $f_6$ & id \\
3.6 & & & id & $f_4$ & id & id & $f_4$ & $f_7$ \\
3.7 & & & & & & & id & id \\
\hline 
\end{tabular}
\end{center}
\label{table:homs}
\end{table}%

Missing edges in the Hasse diagram of $\mathcal{K}_{3,3}$ can be justified by appealing to the following result from \cite{BCDM}.

\begin{prop}\cite{BCDM}\label{prop:Prop}
Let $\overline{G}$ and $ \widehat{G}$ be geometric realizations of a  graph $G$, and suppose $\overline{G}\stackrel{f} \preceq \widehat{G}$. Then each of the following conditions holds.
\begin{enumerate}
  \item\label{1} $\widehat G_0$ is a subgraph of $\G_0$.
  \item\label{2} $f$ induces a graph homomorphism $EX(\G) \to EX(\widehat{G})$;
   \item\label{3} $f$ induces a color-preserving graph homomorphism $LEX(\G) \to LEX(\widehat{G})$  that restricts to an automorphism on $L(G)$.
\end{enumerate}  
\end{prop}

\noindent Part (\ref{1}) and Figure~\ref{fig:allK330} together show that:
\begin{itemize}
\item realization 3.1 does not precede realizations 5.3, 5.4, 5.4, 5.5, 5.6, 5.7 and 5.8; 
\item neither realization 3.2 nor realization 3.3 precedes realizations 5.7 and 5.8; 
\item realizations 5.1, 5.2 and 5.3 do not precede realization 7.2.
\end{itemize}
Part (\ref{2}) and Figure~\ref{fig:lineXK33} together show that:
\begin{itemize}
\item realization 3.5 does not precede realizations 5.1, 5.2 or 5.6; 
\item realization 3.7 does not precede realizations 5.1, 5.2, 5.3, 5.4, 5.5 or 5.6;
\item realizations 5.6, 5.7 and 5.8 do not precede realization 7.1.
\end{itemize} 

\noindent Part (\ref{3}) and Figure~\ref{fig:lineXK33} together show that:
\begin{itemize}
\item realization 3.4 does not precede realization 5.1;
\item realization 3.4 does not precede realization 5.8;
\item realization 3.6 dos not precede realizations 5.1 and 5.2.
\end{itemize}

\bibliographystyle{plain}

\bibliography{Bipartite}

\begin{thebibliography}{1}

\bibitem{BC}
Debra Boutin and Sally Cockburn.
\newblock Geometric graph homomorphisms.
\newblock {\em Journal of Graph Theory}, 69(2):97--113, February 2012.

\bibitem{BCDM}
Debra Boutin, Sally Cockburn, Alice Dean, and Andrei Margea.
\newblock Posets of geometric graphs.
\newblock {\em Ars Mathematica Contemporanea}, 5:265--284, 2012.

\bibitem{HH}
Heiko Harborth.
\newblock Parity of number of crossings for complete $n$-partite graphs.
\newblock {\em Mathematica Slovaca}, 26:77--95, 1976.

\end{thebibliography}

\end{document}